\def\<{\langle}
\def\>{\rangle}

\def \w {\omega}
\def \t {\theta}

\def \b {\beta}

\documentclass[12pt]{iopart}

\usepackage{setstack}
\usepackage{amssymb}
\usepackage{iopams}  
\usepackage{graphicx}
\usepackage{subfigure}
\usepackage{color}

\newtheorem{thm}{Theorem}

\begin{document}
\title{Design of Charge-Balanced Time-Optimal Stimuli for Spiking Neuron Oscillators}
\author{ Isuru S. Dasanayake and Jr-Shin Li}
\address{Department of Electrical \& Systems Engineering,
Washington University in St. Louis, St. Louis, MO 63130, USA.}
\ead{dasanayakei@seas.wustl.edu and jsli@seas.wustl.edu }

\begin{abstract}
	In this paper, we investigate the fundamental limits on how the inter-spike time of a neuron oscillator can be perturbed by the application of a bounded external control input (a current stimulus) with zero net electric charge accumulation. We use phase models to study the dynamics of neurons and derive charge-balanced controls that achieve the minimum and maximum inter-spike times for a given bound on the control amplitude. Our derivation is valid for any arbitrary shape of the phase response curve and for any value of the given control amplitude bound. In addition, we characterize the change in the structures of the charge-balanced time-optimal controls with the allowable control amplitude. 
We demonstrate the applicability of the derived optimal control laws by applying them to mathematically ideal and experimentally observed neuron phase models, including the widely-studied Hodgkin-Huxley phase model, and by verifying them with the corresponding original full state-space models. This work addresses a fundamental problem in the field of neural control and provides a theoretical investigation to the optimal control of oscillatory systems.
	
\end{abstract}

\submitto{\JNE}
\maketitle

\section{Introduction}
Neurons exhibit short-lasting voltage spikes known as action potentials, which are sensitive to external current stimuli \cite{Izhikevich07}. The inter-spike time interval of a neuron characterizes its properties and can be controlled by use of external stimuli. The ability to control neuron spiking activities is fundamental to theoretical neuroscience, and the concept of effective control of such neurological behavior has led to the development of innovative therapeutic procedures \cite{Schiff94,Benabid10} for neurological disorders including deep brain stimulation (DBS) for Parkinson's disease and essential tremor \cite{Lozano04,Nabi11}, where electrical pulses are used to inhibit pathological synchrony among neuron populations. In such neurological treatments and other applications such as the design of artificial cardiac pacemakers \cite{Ortmanns07}, it is of clinical importance to avoid long and strong electrical pulses in order to prevent the tissue from damage,  as well as to maintain zero net electric charge accumulation over each stimulation cycle in order to suppress undesirable side effects. High levels of electric charge accumulation may trigger irreversible electrochemical reactions resulting in damage of neural tissues and corrosion of electrodes \cite{Merrill05}.

Motivated by these practical needs, in this paper we study the design of time-optimal controls for spiking neurons, which lead to the minimum and maximum inter-spike times and remain charge-balanced.
We study the dynamics of neuron oscillators through phase models which are simplified yet accurate models that capture essential overall properties of an oscillating neuron \cite{Izhikevich07, Brown04}, and which form a standard nonlinear control system that characterizes the evolution of an oscillating system by a single variable, namely, the phase. Phase models are conventionally used to investigate the synchronization patterns and study the dynamical responses of oscillators where the inputs to the oscillatory systems are initially defined \cite{Brown04, Ashwin92, Tass89}. Recently, control-theoretic approaches, including calculus of variations and the maximum principle, have been employed to design external stimuli for optimal manipulation of the dynamic behavior of neuron oscillators. These include the design of minimum-power controls for spiking a single neuron at specified time instances \cite{Moehlis06, Dasanayake11, Dasanayake11b}, optimal waveforms for entrainment of neuron ensembles \cite{Zlotnik12,Zlotnik11,kiss10}, and open-loop controls for establishing and maintaining a desired phase configuration, such as anti-phase for two coupled neuron oscillators \cite{Stefanatos12}. 
Work on considering stochastic effects to neuron systems such as the optimal control of neuronal spiking activity receiving a class of random synaptic inputs has also been investigated \cite{Feng03}. In addition, controllability of an ensemble of uncoupled neurons was explored for various mathematically ideal phase models, where an effective computational optimal control method based on pseudospectral approximations was employed to construct optimal controls that elicit simultaneous spikes of a neuron ensemble \cite{Li_NOLCOS10, Li12}. The derivation of time-optimal and spike timing controls for spiking neurons has been attempted for limited classes of control functions \cite{Nabi10, Danzl08}, however, a complete characterization of the optimal solutions has not been provided, and an analytical and systematic approach for synthesizing the time-optimal controls has been missing. 

In this paper, we derive charge-balanced time-optimal controls for a given bound on the control amplitude and fully characterize the possible range of neuron spiking times determined by such optimal controls. Employing techniques from the optimal control theory, we are able to reveal different structures of the time-optimal controls that vary with the allowable bound of the control amplitude. Moreover, we validate these controls derived according to  phase models by applying them to the corresponding original full state-space neuron models. As a demonstration, the validation is performed using the Hodgkin-Huxley equations \cite{Hodgkin52}, where the spiking behavior of the state-space model shows great qualitative agreement with that of the phase model and which demonstrates the applicability of our theoretical results based on the phase model. Such an important validation, which is largely lacking in the literature, allows us to explore the fundamental limits of the phase model as an approximation of state-space models.

This paper is organized as follows. In Section \ref{sec:toc}, we consider the time-optimal control of a general phase oscillator and derive the charge-balanced minimum-time and maximum-time controls with constrained control amplitude by using the Pontryagin's maximum principle \cite{Pontryagin62}. In Section \ref{sec:example}, we apply the derived optimal control strategies to both mathematically ideal and experimentally observed phase models, including the well-known SNIPER \cite{Brown04}, Hodgkin-Huxley, and Morris-Lecar \cite{Morris81} phase models, and present the simulated optimal solutions. In Section \ref{sec:validation}, we validate the obtained optimal controls through the Hodgkin-Huxley model.

\section{Charge-Balanced Time-Optimal Control for Phase Models of Spiking Neuron Oscillators}
\label{sec:toc}
The dynamics of a periodically spiking neuron oscillator can be described by a phase model of the form \cite{Brown04}
\begin{equation}
    \label{eq:phasemodel}
    \frac{d\theta}{dt}=\omega+Z(\theta)u(t),
\end{equation}
where $\theta$ denotes the phase of the oscillation, $\omega>0$ is neuron's natural oscillation frequency, and $u(t)\in\mathbb{R}$ is the external current stimulus (control) that is applied to perturb the phase dynamics of the neuron. The real-valued function $Z(\theta)$ is the phase response curve (PRC) that characterizes the infinitesimal change of the phase to an external control input. Conventionally, the neuron is said to spike when its phase $\theta=2n\pi$, where $n\in\mathbb{N}$. In the absence of any input $u(t)$, the neuron spikes periodically at its natural frequency, while the spiking time may be advanced or delayed in a desired manner by the application of an appropriate weak control.

\subsection{Charge-Balanced Minimum-Time Control}
The optimal design of controls that lead to the minimum inter-spiking time of a neuron subject to a given bound on the control amplitude and the charge-balance constraint can be formulated as a time-optimal steering problem of the form 
\begin{eqnarray}
    \min_{u(t)} \quad & T, \nonumber\\
    {\rm s.t.} \quad & \dot{\theta}=\omega+Z(\theta)u(t), \nonumber\\
		\label{eq:toc}
    &\dot{p}=u(t),\\
    &|u(t)|\leq M, \ \ \forall\ t\in[0,T], \nonumber\\
    &\theta(0)=0, \quad \theta(T)=2\pi, \nonumber\\
    &p(0)=0, \quad p(T)=0, \nonumber
\end{eqnarray}
where $T$ is the inter-spiking time required that we wish to minimize and $M>0$ is the bound of the control amplitude. The constraints involving the time-dependent variable $p(t)$ are equivalent to the charge-balance constraint, i.e., $p(t)=\int_0^t u(\sigma)d\sigma=0$ with $p(0)=p(T)=0$, guaranteeing that the charge accumulated over a spiking cycle is zero. Note that the consideration of bounded controls is of fundamental importance since the phase reduction is no longer valid when the control exceeds a level that can be considered weak.


\subsubsection{Derivation of the Charge-Balanced Minimum-Time Control:} \label{sec:deriv_min}
The Hamiltonian of the optimal control problem as in (\ref{eq:toc}) is given by
\begin{equation}
	\label{eq:hamiltonian}
	H=\lambda_0+\lambda_1(\omega+Z(\theta)u)+\lambda_2 u
\end{equation}
where $\lambda_0\geq 0,\ \lambda_1$, and $\lambda_2$ are Lagrange multipliers associated with the Lagrangian, system dynamics, and the charge-balance constraint, respectively. According to the optimality conditions of the maximum principle (see \ref{Appendix1}), the adjoint variables $\lambda_1$ and $\lambda_2$ must satisfy the time-varying equations  $\dot{\lambda}_1=-\frac{\partial H}{\partial \theta}$ and $\dot{\lambda}_2=-\frac{\partial H}{\partial p}$, respectively, which yields
\begin{eqnarray}
	\dot{\lambda_1}=&-\lambda_1 u \frac{\partial Z(\theta)}{\partial\theta}, \label{eq:adj1}\\
	\dot{\lambda_2}=&0, \label{eq:adj2}
\end{eqnarray}
and hence $\lambda_2$ is a constant. Since the Hamiltonian $H$ is not explicitly dependent on time and the terminal time is free, we have $H\equiv0$, $\forall\, t\in[0,T]$, along the optimal trajectory from the maximum principle.

It is straightforward to see from a rearrangement of \eref{eq:hamiltonian}, $H=\lambda_0+\lambda_1\w+(\lambda_1 Z(\t)+\lambda_2)u$, that the control
\begin{equation}
	\label{eq:I*Tmin_bang}
    u^*_{min}=\left\{\begin{array}{ll} M, & \phi< 0 \\ -M, & \phi\geq0
    \end{array}\right.
\end{equation}
minimizes the Hamiltonian $H$, where 
\begin{eqnarray}
	\label{eq:phi}
	\phi=\lambda_1 Z+\lambda_2
\end{eqnarray}
is called the switching function. Hence, according to the maximum principle, $u^*_{min}$ is a candidate of the optimal solution to the problem as in \eref{eq:toc}, provided $\phi=0$ for a nonzero time period does not occur. This type of controls is known as bang-bang controls, which takes only the extremals of the control set, e.g., $-M$ or $M$ in this case. The switching between $-M$ and $M$ occurs at $\phi=0$ and the challenge is to calculate the values of the  multipliers $\lambda_1$ and $\lambda_2$, which define the function $\phi$ and thus the optimal control sequence.

An alternative candidate of the minimum-time control may exist. If $\phi\equiv 0$ for some non-zero time interval $S=[\tau_1,\tau_2]$, then its derivatives $\dot{\phi}$, $\ddot{\phi}$, etc., will also be equal to zero over $S$. In this case, the bang-bang control \eref{eq:I*Tmin_bang} may not be optimal. Such a control that forces the switching function $\phi$ and all of its derivatives to vanish over a time period is known as a singular control \cite{Bonnard03}, and it can be calculated according to the following fashion. When $\phi=0$, $\dot{\phi}=0$, $\ddot{\phi}=0$, $\ldots$, for a given time interval $S$, we have
\begin{equation}
	\label{eq:phi0}
	\phi=\lambda_1 Z+\lambda_2=0
\end{equation}
and then, by substituting from \eref{eq:phasemodel}, \eref{eq:adj1}, and \eref{eq:adj2}, the function $\dot{\phi}$ is given by
\begin{equation}
	\label{eq:dotphi0}
	\dot{\phi}=\lambda_1\omega \frac{\partial Z}{\partial \theta}=0
\end{equation}
which yields $\frac{\partial Z}{\partial \theta}=0$ because $\omega>0$ and $\lambda_1\neq0$. The latter is due to the non-triviality condition of the maximum principle, i.e., $(\lambda_0,\lambda_1,\lambda_2)\neq\mathbf{0}$, since $\lambda_2=0$ if $\lambda_1=0$ from \eref{eq:phi0}, which leads to $\lambda_0=0$ from \eref{eq:hamiltonian} as $H\equiv0$. Therefore, $\lambda_1\neq0$ holds along the optimal trajectory and $\frac{\partial Z}{\partial \theta}=0$ defines a singular trajectory, i.e., the trajectory of the system following a singular control. As in the calculation of $\dot{\phi}$, the second derivative $\ddot{\phi}$ can be obtained using \eref{eq:phasemodel} and $\frac{\partial Z}{\partial \theta}=0$ to get
\begin{equation}
	\label{eq:dotdotphi0}
	\ddot{\phi}=\lambda_1\omega \frac{\partial^2 Z}{\partial \theta^2}(\omega+Zu).
\end{equation}
It is clear from \eref{eq:dotdotphi0} that if $ \frac{\partial^2 Z}{\partial \theta^2}\neq 0$, the control that makes $\ddot{\phi}=0$ is given by $u_{s}=-\omega/Z$. In the case when $\frac{\partial^2 Z}{\partial \theta^2}=0$, we need to calculate $\dddot{\phi}$ in order to get the singular control $u_{s}$. However,  no matter how many derivatives are used, the singular control is given by the same form, $u_{s}=-\omega/Z$.

If a singular trajectory exists, then one must examine whether it is ``fast" or ``slow" compared to the bang-bang trajectory in order to determine the minimum-time control. Suppose that the singular control $u_{s}=-\omega/Z$ is admissible over a nonzero time interval $S=[\tau_1,\tau_2]$. Then, from \eref{eq:phasemodel} the phase velocity is equal to zero, i.e., $\dot{\theta}\equiv 0$, over $S$ by the application of $u_{s}$. This implies that the singular trajectory is slower than any feasible trajectory along which $\dot{\theta}\geq0$ over $S$. 
Therefore, the charge-balanced control that spikes neurons in minimum time is of the bang-bang form.


\subsubsection{Computation and Synthesis of the Charge-Balanced Minimum-Time Control:} \label{sec:comp_min}
Because the minimum spiking time of the neuron system as in \eref{eq:phasemodel} is achieved by a bang-bang control, it remains to calculate the switching points in order to synthesize this time-optimal control. Since $\phi=0$ holds at the switching points, according to \eref{eq:phi0}, these points are defined via the inverse function of the PRC,
\begin{equation}
	\label{eq:switch}
	\theta_s=Z^{-1}\left(-\frac{\lambda_2}{\lambda_1}\right).
\end{equation}
In addition, with the Hamiltonian condition $H\equiv0$, the value of the multiplier $\lambda_1$ is given by $\lambda_1=-\frac{\lambda_0}{\omega}$ at these switching points. Without loss of generality, we let $\lambda_0=1$, which leads to $\lambda_1=-\frac{1}{\omega}$. Applying this to \eref{eq:switch} results in
\begin{eqnarray}
	\label{eq:switch_modified}
	\theta_s=Z^{-1}(\alpha)=Z^{-1}\left(\lambda_2 \omega\right),
\end{eqnarray}
where $\lambda_2$ and $\omega$ are both constants. Let $Z^{-1}(\alpha)$ have $n$ solutions in the interval $(0,2\pi)$ given by $\theta_1,\theta_2,\ldots\theta_n$, and define $\theta_0=0$ and $\theta_{n+1}=2\pi$. Then, if we start with the control $u=M$, the charge-balance constraint gives rise to the condition
\begin{equation}
	\label{eq:charge}
	0=\int_0^T u(t)dt=\sum_{i=0}^{i=n}\int_{\theta_i}^{\theta_{i+1}}{\frac{(-1)^iM}{\omega+(-1)^iZ(\theta)M}}d\theta
\end{equation}
and the total time $T$ under this bang-bang control is represented by
\begin{equation}
	\label{eq:time}
	T=\sum_{i=0}^{i=n}\int_{\theta_i}^{\theta_{i+1}}{\frac{1}{\omega+(-1)^iZ(\theta)M}}d\theta.
\end{equation}
Equation \eref{eq:charge} together with the switching conditions $Z(\theta_i)=\alpha$ for $i=1,2,\ldots n$ define $n+1$ equations of $n+1$ variables, $\left\{\theta_1,\theta_2,\ldots\theta_n,\alpha\right\}$. This system of equations can be solved to get the set of optimal switching angles, denoted as $S_M$, and the constant $\alpha$. Similarly, if we start  with the control $u=-M$, by substituting $M$ with $-M$ in \eref{eq:charge} we obtain the other set of solutions, denoted as $S_{-M}$. The bang-bang control, determined by the set of switching angles, which results in the shorter spiking time is the charge-balanced minimum-time control, while the opposite case is a candidate for the charge-balanced maximum-time control.

Alternatively, given the two sets of switching angles, the optimal switching sequence can be determined by computing $\dot{\phi}$ at the switching points. We denote the vector fields corresponding to the constant bang controls $u(t)\equiv -M$ and $u(t)\equiv M$ by $X = \omega - MZ$ and $Y = \omega + MZ$, respectively, and call the respective trajectories corresponding to them as $X$- and $Y$- trajectories. A concatenation of an $X$-trajectory followed by a $Y$-trajectory is denoted by $XY$, while the concatenation in the reverse order is denoted by $YX$. If $\dot{\phi}<0$ at a switching point, then the $X$ to $Y$ switch is optimal according to the switching law \eref{eq:I*Tmin_bang}, and similarly if $\dot{\phi}>0$, then the $Y$ to $X$ switch is optimal. Since $\lambda_1=-1/\omega$ at the switching points, we have
\begin{equation}
	\label{eq:dot_phi_at_switch}
	\dot{\phi}=\lambda_1\omega \frac{\partial Z}{\partial \theta}=-\frac{\partial Z}{\partial \theta}.
\end{equation} 
Therefore, the value of $\frac{\partial Z}{\partial \theta}$ at the switching points defines the switching type. If $\frac{\partial Z}{\partial \theta}>0$, an $X$ to $Y$ switch is optimal and if $\frac{\partial Z}{\partial \theta}<0$, a $Y$ to $X$ switch is optimal. 

\subsection{Charge-Balanced Maximum-Time Control}
\label{sec:deriv_max}
\subsubsection{(CaseI: Bang-Bang Control)\label{sec:deriv_max_bb}} When the control amplitude is limited by $M<\min\big\{\big|\frac{\omega}{Z(\theta)}\big|:\theta\in[0,2\pi)\big\}$, singular controls are not admissible since $u_{s}=-\omega/Z$ as shown in Section \ref{sec:deriv_min}. Therefore, the maximum-time control is given by the bang-bang form
\begin{equation}
    u^*_{max}=\left\{\begin{array}{ll} -M, & \phi\leq 0 \\ M, & \phi>0.\end{array}\right.
\end{equation}
where $\phi$ is defined as in \eref{eq:phi}.  
The optimal switching sequence is determined between $S_M$ and $S_{-M}$, whichever results in longer spiking time. Another way to determine the optimal switching sequence is by evaluating $\frac{\partial Z}{\partial \theta}$ at the switching points as described in Section \ref{sec:comp_min}. When $\frac{\partial Z}{\partial \theta}>0$ at a switching point, a $Y$ to $X$ switch is optimal,  while when $\frac{\partial Z}{\partial \theta}<0$, an $X$ to $Y$ switch is optimal. 

\subsubsection{(CaseII: Bang-Singular-Bang Control)\label{sec:deriv_max_bsb}}  When singular controls are admissible, that is, when $M\geq\min\big\{\big|\frac{\omega}{Z(\theta)}\big|:\theta\in[0,2\pi)\big\}$, the maximum-time control is a combination of bang and singular controls (see the examples in Section \ref{sec:max_time_sniper} and \ref{sec:HH}). The procedure of the optimal control synthesis is to choose a bang control that drives the system to a singular trajectory (a system trajectory following a singular control), staying on that trajectory, and then exiting at the point from which a bang control can steer the system to the desired terminal state. Examples involving the construction of charge-balanced minimum-time and maximum-time optimal controls are illustrated in Section \ref{sec:example}.

\section{Examples}
\label{sec:example}
We now apply the derived optimal control strategies to several commonly-used phase models characterized by various PRC's, including mathematically ideal and experimentally observed phase models. These examples demonstrate the applicability of our optimal control methods to manipulate neuron dynamics. We emphasize that these optimal controls are designed with respect to a given bound of the control amplitude, so that they can be made practical and satisfy the weak forcing assumption in the phase model.

\subsection{SNIPER Phase Model}
\label{sec:SNIPER}
The SNIPER phase model is characterized by a type I PRC and is of the form \cite{Brown04}
\begin{eqnarray}
	\label{eq:sniper_phasemodel}
	\dot{\theta}=\omega+z_d\, (1-\cos\theta) u,
\end{eqnarray}
where $\omega$ is the natural oscillation frequency of the neuron, $z_d$ is a model-dependent constant, and $u$ is the external stimulus. This model is derived from a SNIPER bifurcation (saddle-node bifurcation of a fixed point on a periodic orbit), which can be found on Type I neurons \cite{Ermentrout96} such as the Hindmarsh-Rose model \cite{Rose89}. Neurons described by this model spike periodically with the natural period $T_0=2\pi/\omega$ in the absence of any external input $u$. 

Before calculating the minimum- and maximum-time spiking controls for the SNIPER phase model, we first examine the existence of singular trajectories. Recall from \eref{eq:dotphi0} that the singular trajectory is defined by $\frac{\partial Z}{\partial \theta}=0$, which yields
$$z_d\sin\theta=0.$$
Therefore, there exist three possible singular trajectories (in this case singular points), $\theta=0$, $\theta=\pi$, and $\theta=2\pi$. The points $\theta=0$ and $\theta=2\pi$ are infeasible singular points, at which the nonzero phase velocity, $\dot{\theta}=\omega$, immediately forces the system away from these points, Hence, $\theta_s=\pi$ is the only possible singular point, and the singular control $u=-\omega/Z(\theta_s)=-\omega/(2z_d)$ yields $\dot{\theta}=0$ at $\theta_s$, making the system stay at $\theta_s$.


\subsubsection{Charge-Balanced Minimum-Time Control for SNIPER Phase Model:}
\label{sec:min_time_sniper}
\begin{figure}
	\centering
		\includegraphics[scale=0.5]{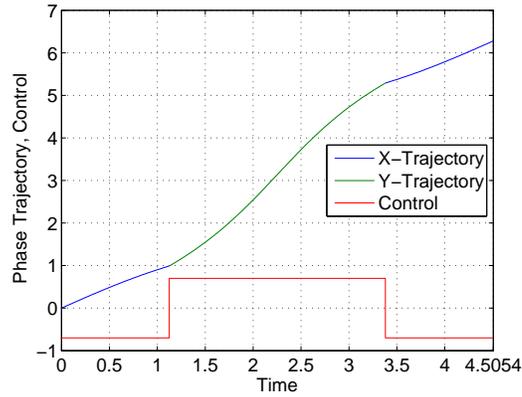}
	\caption{The charge-balanced minimum-time control and the corresponding phase trajectory for the SNIPER phase model with $z_d=1$, $\omega=1$, and $M=0.7$.}
	\label{fig:sniper_min_time}
\end{figure}

Since the charge-balanced minimum-time control takes the bang-bang form as shown in Section \ref{sec:deriv_min}, the switching points are given from \eref{eq:switch_modified} by
\begin{equation}
	\theta_s=\cos^{-1}\left\{1-\frac{\omega\lambda_2}{z_d}\right\}.
\end{equation}
The cosine function has two solutions in $[0,2\pi)$ and thus there are two switching points $\theta_1=\gamma$ and $\theta_2=2\pi-\gamma$ with $\gamma\in[0,\pi)$. Because $\lambda_1=-1/\omega$ for both switching points and the derivative of the switching function $\dot{\phi}=-z_d\sin\theta<0$ for $\theta\in(0,\pi)$, if a switch occurs on the interval $(0,\pi)$, it will be $X$ to $Y$. Reversely, if a switch occurs on $(\pi,2\pi)$, then it will be $Y$ to $X$ because $\dot{\phi}>0$ for $\theta\in(\pi,2\pi)$. It follows that an $XYX$ trajectory is optimal for achieving the minimum inter-spike time. The parameter $\gamma$ that defines the switching points is calculated using the charge-balance constraint as in \eref{eq:charge} by solving $R(M,\gamma)=0$, where
\begin{equation}
	\label{eq:charge_sniper}
	R(M,\gamma)=\int_{0}^{\gamma}{\frac{-M}{\omega-z_d(1-\cos\theta)M}}d\theta+\int_{\gamma}^{\pi}{\frac{M}{\omega+z_d(1-\cos\theta)M}}d\theta.
\end{equation}
 Then, the optimal control is given by
\begin{equation}
\label{eq:opt_control_sniper}
    u^*_{min}=\left\{\begin{array}{ll} -M, & 0\leq \theta< \theta_1 ,\\ M, & \theta_1\leq\theta\leq\theta_2, \\ -M, & \theta_2<\theta\leq2\pi,
    \end{array}\right.
\end{equation}
 and by following \eref{eq:time} the time required to spike the neuron, namely, to reach $\theta=2\pi$, is given by
\begin{equation}
\label{eq:opt_time_sniper}
    T=\int_{0}^{\gamma}{\frac{4}{\omega-z_d(1-\cos\theta)M}}d\theta.
\end{equation}
\Fref{fig:sniper_min_time} shows the charge-balanced minimum-time control and the corresponding phase trajectory for the SNIPER phase model with $z_d=1$, $\omega=1$, and $M=0.7$.

\subsubsection{Charge-Balanced Maximum-Time Control for SNIPER Phase Model:}
\label{sec:max_time_sniper}
\begin{figure}
	\centering
		\includegraphics[scale=0.5]{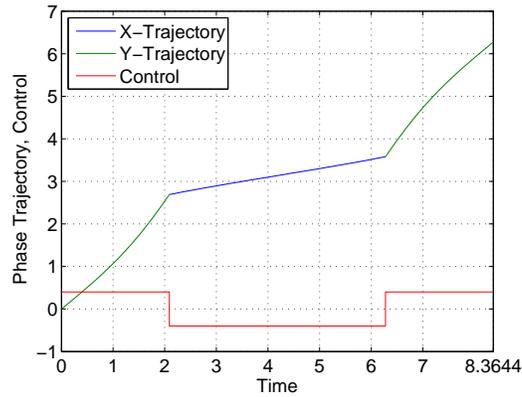}
	\caption{The charge-balanced maximum-time control and the corresponding phase trajectory for the SNIPER phase model with $z_d=1$, $\omega=1$, and $M=0.4<\frac{\omega}{2z_d}=0.5$.}
	\label{fig:sniper_max_time_bang_bang}
\end{figure}

 There are two control scenarios for maximizing the spiking time of a SNIPER neuron depending on the control amplitude.
\ \\
\noindent(\emph{Case I: $M<\frac{\omega}{2z_d}$}) If the bound of the control amplitude $M<|\frac{\omega}{Z(\theta)}|=|\frac{\omega}{z_d(1-\cos\theta)}|<\frac{\omega}{2z_d}$, then there exist no admissible singular controls and the maximum-time control takes the bang-bang form as described in Section \ref{sec:deriv_max_bb}. In this case, there are two switches and the $YXY$ trajectory is optimal. The maximum-time control is given by
\begin{equation}
\label{eq:max_time_control_sniper}
    u^*_{max}=\left\{\begin{array}{ll} M, & 0\leq\theta< \theta_1 ,\\ -M, & \theta_1\leq\theta\leq\theta_2, \\ M, & \theta_2<\theta\leq2\pi,
    \end{array}\right.
\end{equation}
where $\theta_1=\beta$, $\theta_2=2\pi-\beta$, and the parameter $\beta$ is obtained by solving $R(-M,\b)$ as defined in \eref{eq:charge_sniper}.
\Fref{fig:sniper_max_time_bang_bang} illustrates the charge-balanced maximum-time control and the corresponding phase trajectory for the SNIPER phase model with $z_d=1$, $\omega=1$, and $M=0.4<\frac{\omega}{2z_d}=0.5$. \\

\begin{figure}
	\centering
		\includegraphics[scale=0.5]{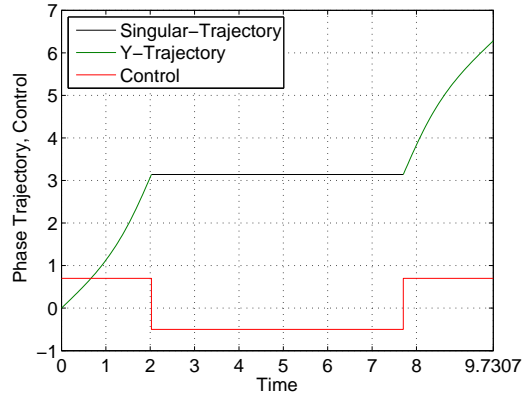}
	\caption{ The maximum-time charge-balanced control and corresponding phase trajectory for the SNIPER phase model with $z_d=1,\omega=1$, and $M=0.7$.}
	\label{fig:sniper_max_time_bang_sing_bang}
\end{figure}

\noindent(\emph{Case II: $M\geq\frac{\omega}{2z_d}$}) In this case, the system can be driven along the singular trajectory which is optimal (slower than the bang control), and the maximum-time control takes the bang-singular-bang form. Because, for example, when $\theta\in(0,\pi)$, the derivative of the switching function $\dot{\phi}=-z_d\sin\theta<0$, and then the $YX$ trajectory is a candidate for optimality if a switch occurs. However, following an $X$-trajectory with $u=-M\leq\frac{-\omega}{2z_d}$, the singular point $\theta=\pi$ is unreachable. Hence, switching in the interval $(0,\pi)$ is not allowed, and the $Y$-trajectory is optimal for $\theta\in[0,\pi)$. The same reasoning applies for the regime $\theta\in(\pi,2\pi]$, where $Y$-trajectory again is optimal. As a result, the optimal control is of the ``$Y$-singular-$Y$'' form given by
\begin{equation}
\label{eq:max_time_control_sniper_singular}
    u^*_{max}=\left\{\begin{array}{ll} M, & 0\leq\theta< \pi ,\\ -\frac{\omega}{2z_d}, & \theta=\pi, \\ M, & \pi<\theta\leq2\pi.
    \end{array}\right.
\end{equation}
Because $\dot{\theta}=0$ holds along the singular trajectory (in this case the singular point $\theta_s=\pi$), the time duration over which the system stays on it is calculated according to the charge-balance constraint. Let $t_1$ and $t_2$ denote the times for which the first bang control and the singular control are applied, respectively. Since $t_1$ is the time that the system takes to reach $\theta_s=\pi$ by a $Y$-trajectory, we have
\begin{equation}
	t_1=\int_0^\pi{\frac{1}{\omega+z_d(1-\cos\theta)M}}d\theta.
\end{equation}
By symmetry, the amount of time that the system takes following a $Y$-trajectory from $\theta=\pi$ to $\theta=2\pi$ is also $t_1$. Then, $t_2$ is given by
$$t_2=\frac{4Mz_dt_1}{\omega}$$
in order to fulfill the charge-balance constraint. Now the charge-balanced maximum-time control can be stated in terms of time as
\begin{equation}
	\label{eq:max_time_control_sniper_singular_vs_time}
	u^*_{max}=\left\{\begin{array}{ll} M, & 0\leq t< t_1 ,\\ -\frac{\omega}{2z_d}, & t_1\leq t\leq t_1+t_2, \\ M, & t_1+t_2<t\leq t_2+2t_1.
\end{array}\right.
\end{equation} 
\Fref{fig:sniper_max_time_bang_sing_bang} shows the maximum-time charge-balanced control and the corresponding phase trajectory for the SNIPER phase model with $z_d=1$, $\omega=1$, and $M=0.7\geq\frac{\omega}{2z_d}=0.5$.

In the following, we demonstrate the robustness of our analytical method to construct optimal controls for spiking neurons of arbitrary practical PRCs through the Hodgkin-Huxley and Morris-Lecar phase models.

\subsection{Hodgkin-Huxley Phase Model}
\label{sec:HH}
\begin{figure}
	\centering
	\begin{tabular}{cc}
		\subfigure[]{\includegraphics[scale=0.5]{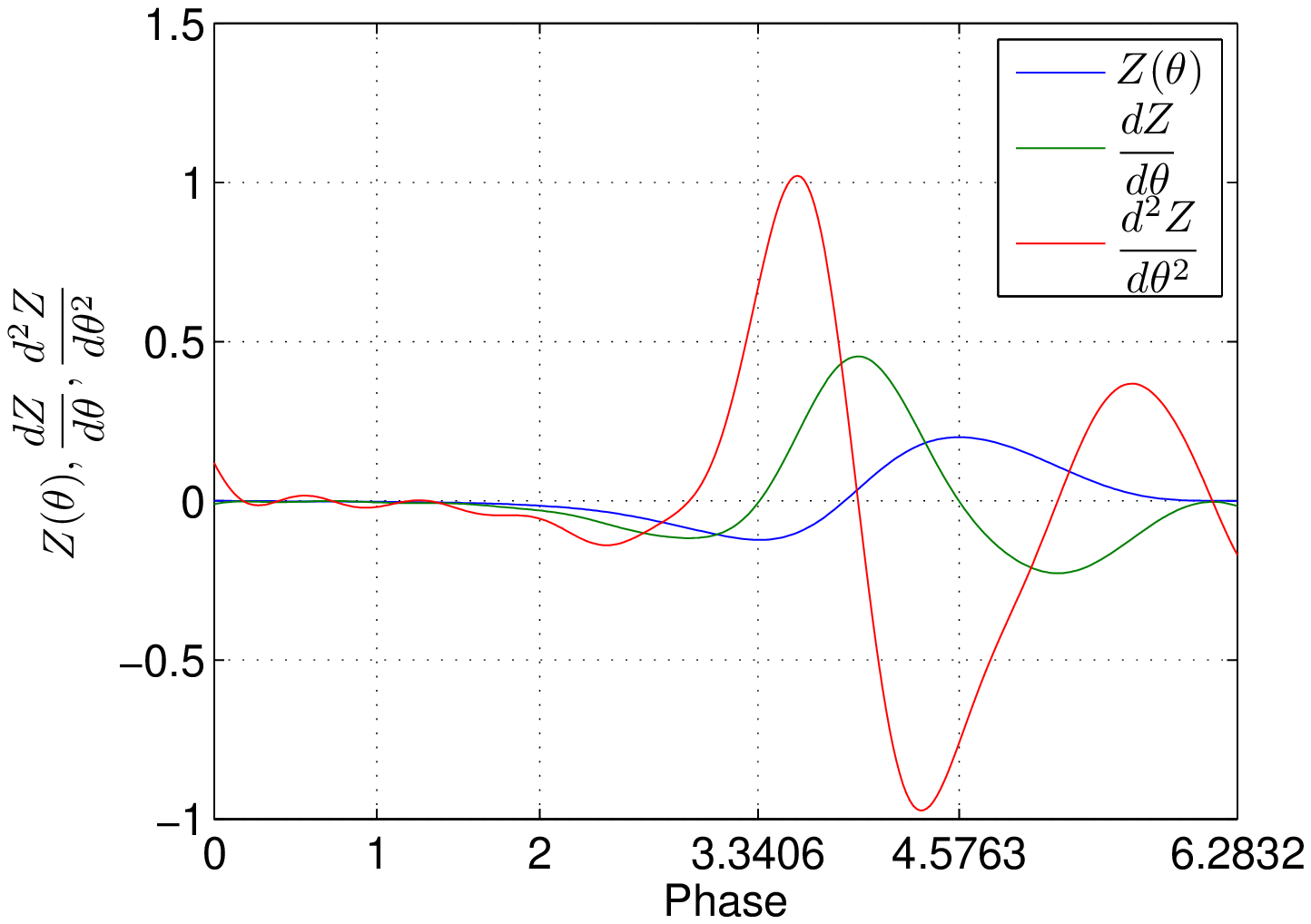}\label{fig:HH_PRC}} &
		\subfigure[]{\includegraphics[scale=0.5]{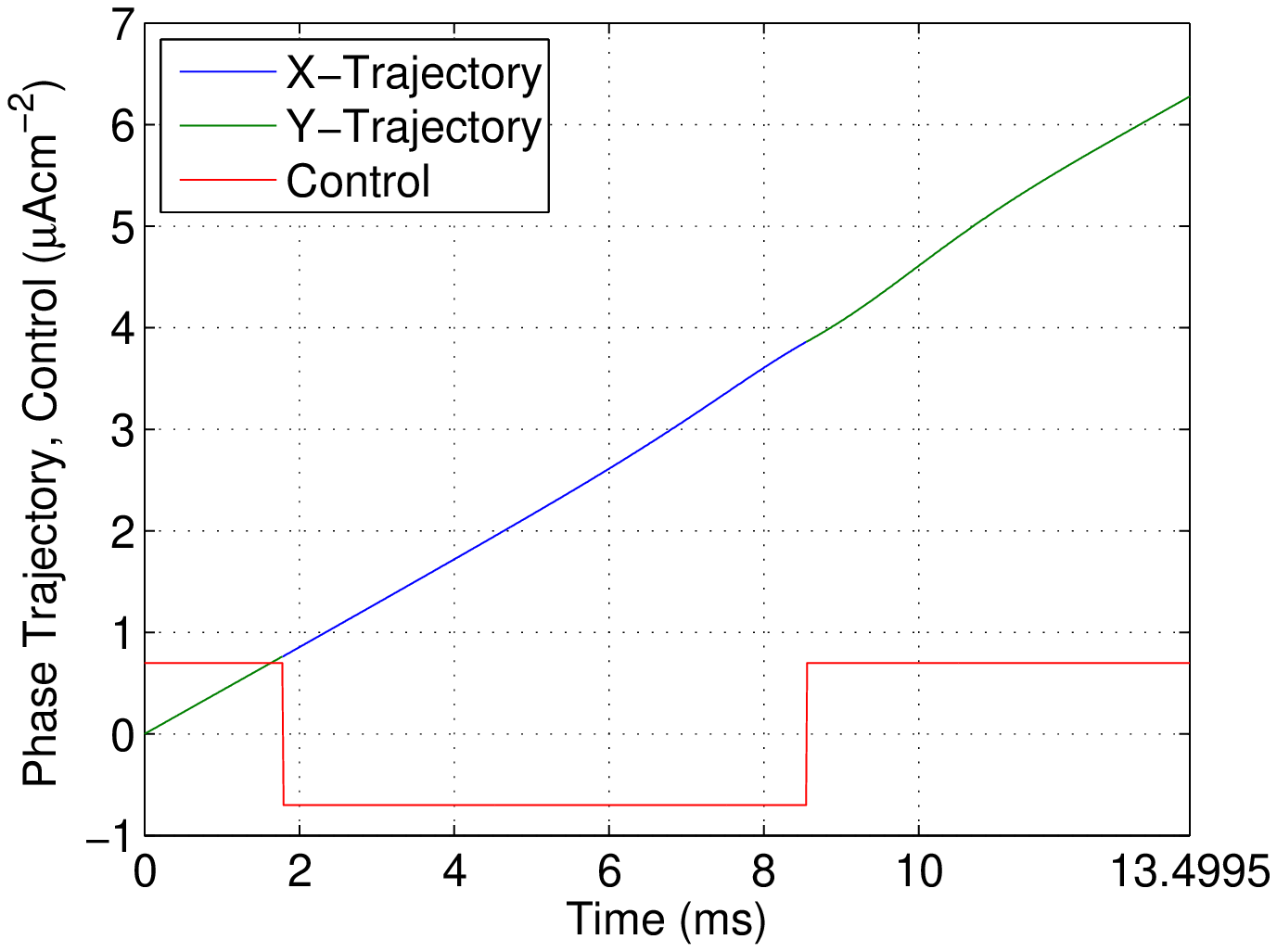}\label{fig:HH_min_time}} \\
		\subfigure[]{\includegraphics[scale=0.5]{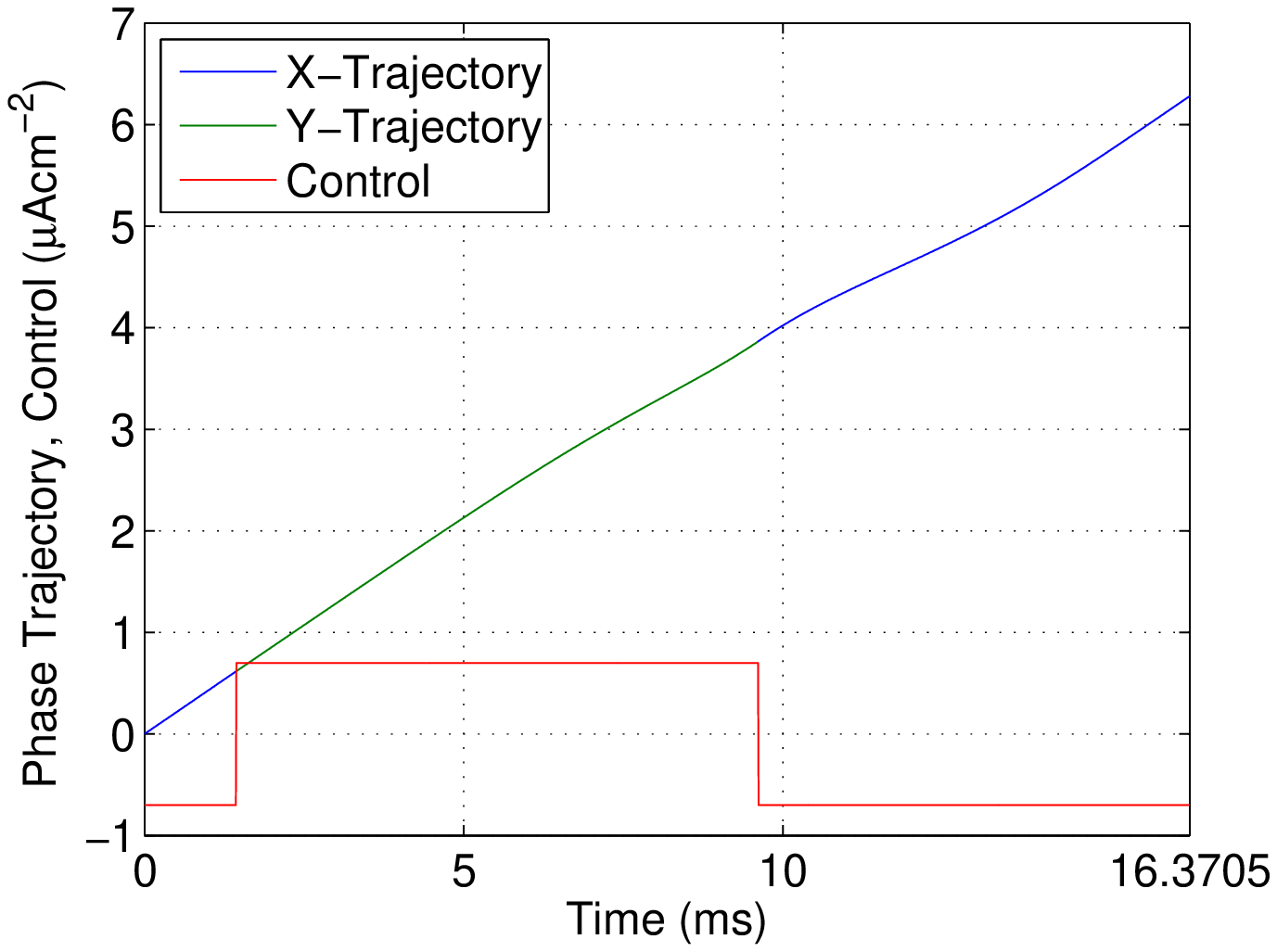} \label{fig:HH_max_time_bang_bang}} &
		\subfigure[]{\includegraphics[scale=0.5]{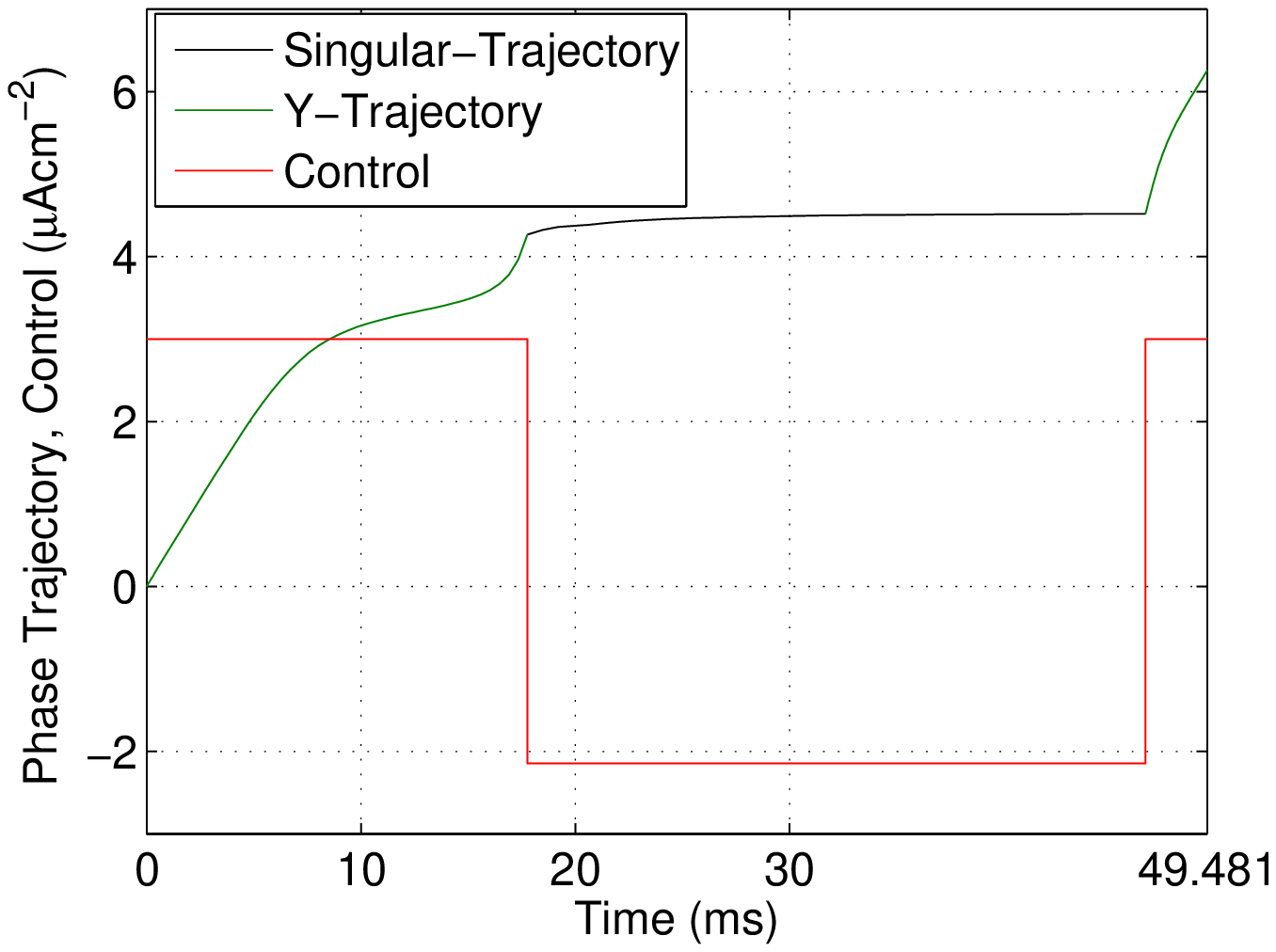} \label{fig:HH_max_time_bang_sing_bang}}
	\end{tabular}
	\caption{\ref{fig:HH_PRC} The Hodgkin-Huxley PRC $Z(\theta)$ and its derivatives, $\frac{dZ}{d\theta}$ and $\frac{d^2Z}{d\theta^2}$. \ref{fig:HH_min_time} The charge-balanced minimum-time control and the corresponding phase trajectory for the Hodgkin-Huxley phase model with respect to the bound on the control amplitude $M=0.7\ \mu Acm^{-2}$.  \ref{fig:HH_max_time_bang_bang} and \ref{fig:HH_max_time_bang_sing_bang} show the charge-balanced maximum-time controls and the corresponding phase trajectories for $M=0.7\ \mu A cm^{-2}$ and $M=3.0\ \mu A cm^{-2}$.} 
\end{figure}

The Hodgkin-Huxley model is a nonlinear system that characterizes the propagation and initiation of the action potential in a squid axon \cite{Hodgkin52}. For the set of parameter values given in \cite{Brown04}, the system exhibits periodic motion with natural frequency $\omega=0.43\, rad/ms$. Its PRC and the first and second derivatives of the PRC are depicted in \Fref{fig:HH_PRC}. To proceed the calculations, we approximate the numerically obtained PRC with eight harmonic terms given by
\begin{eqnarray}
	\label{eq:PRC_coeff}
	Z(\theta)=\sum_{i=1}^{8}{a_i\sin(b_i\theta+c_i)},
\end{eqnarray}
where the coefficients $a_i,b_i$ and $c_i$ are obtained by least squares fit and given in \Tref{Table:HH_coefficient}. In this case, there are two possible singular points, $\theta_{s,1}=3.34$ and $\theta_{s,2}=4.58$, satisfying $\partial Z(\theta)/\partial\theta=0$.

\begin{table}
		\centering
		\begin{tabular}{|c|c|c|c|c|c|c|c|c|}
		\hline
		i & 1 & 2 & 3 & 4 & 5 & 6 & 7 & 8\\
		\hline
			$a_i$ & 0.09176 & 0.07462 & 0.03807 & 0.02425 & 0.01747 & 0.006474 & 0.002752 & 0.0008111 \\
			$b_i$ & 1.002 & 1.996 & 3.002 & 0.5 & 3.747 & 3.747 & 6.228 & 7.651 \\
			$c_i$ & 2.609 & -1.605 & 0.7233 & 0.5148 & 3.552 & -0.7648 & 0.6429 & -4.726 \\
			\hline
		\end{tabular}
		\caption{The coefficients of the equation \eref{eq:PRC_coeff} for the Hodgkin-Huxley PRC.}
		\label{Table:HH_coefficient}
\end{table}
The charge-balanced minimum-time control, which is of the $YXY$ form, and the resulting phase trajectory for the control amplitude bound $M=0.7 \mu Acm^{-2}$ are shown in \Fref{fig:HH_min_time}. The charge-balanced maximum-time controls can take the bang-bang or the bang-sigular-bang form depending on the bound $M$. The cases for $M=0.7 \mu Acm^{-2}$ and $M=3.0 \mu Acm^{-2}$ are illustrated in \Fref{fig:HH_max_time_bang_bang} and \Fref{fig:HH_max_time_bang_sing_bang}, respectively. The  detailed derivations of these optimal controls are presented in \ref{Appendix_HH}.

\subsection{Morris-Lecar Phase Model}
\begin{figure}
	\centering
	\begin{tabular}{cc}
		\subfigure[]{\includegraphics[scale=0.5]{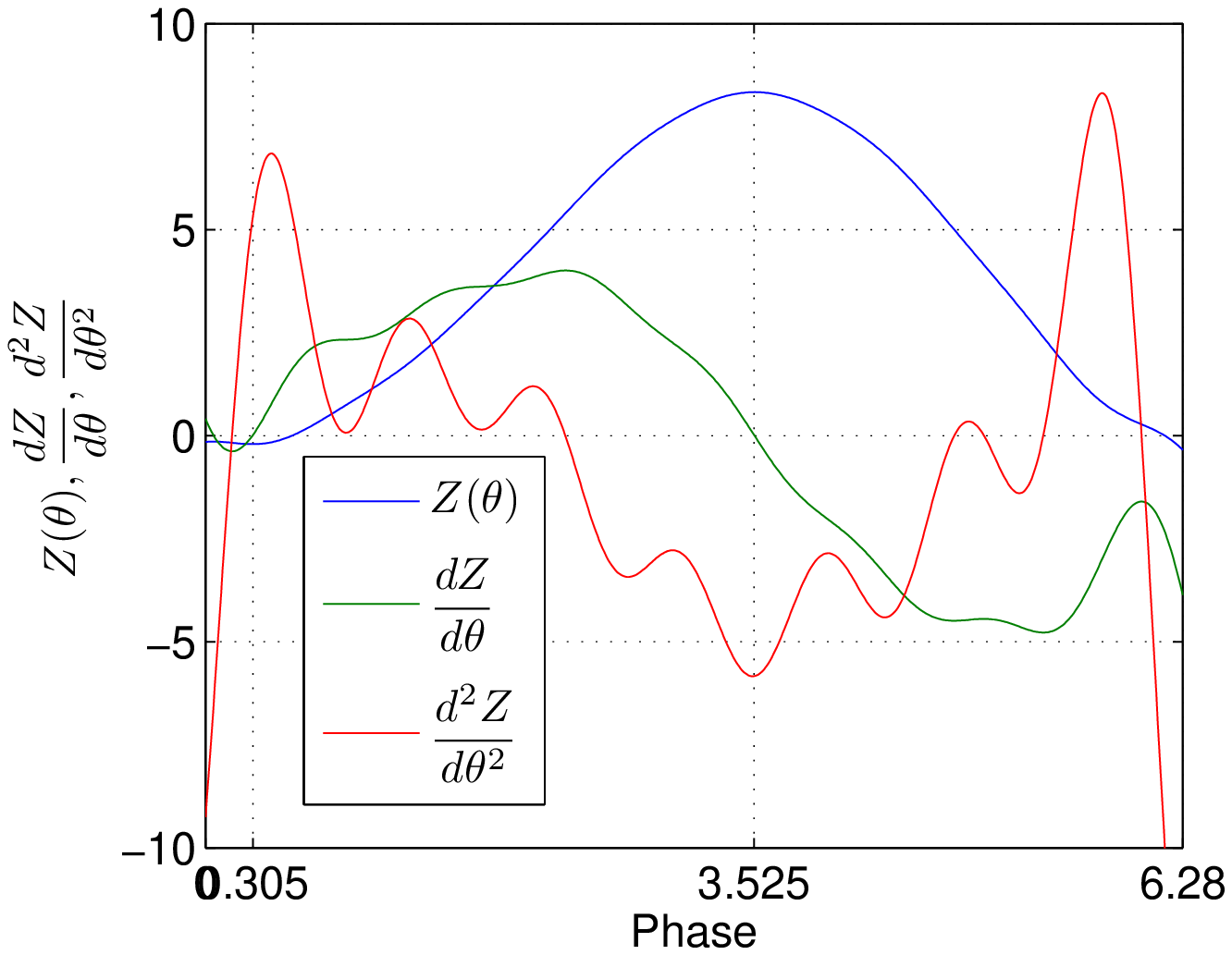}\label{fig:ML_PRC}} &
		\subfigure[]{\includegraphics[scale=0.5]{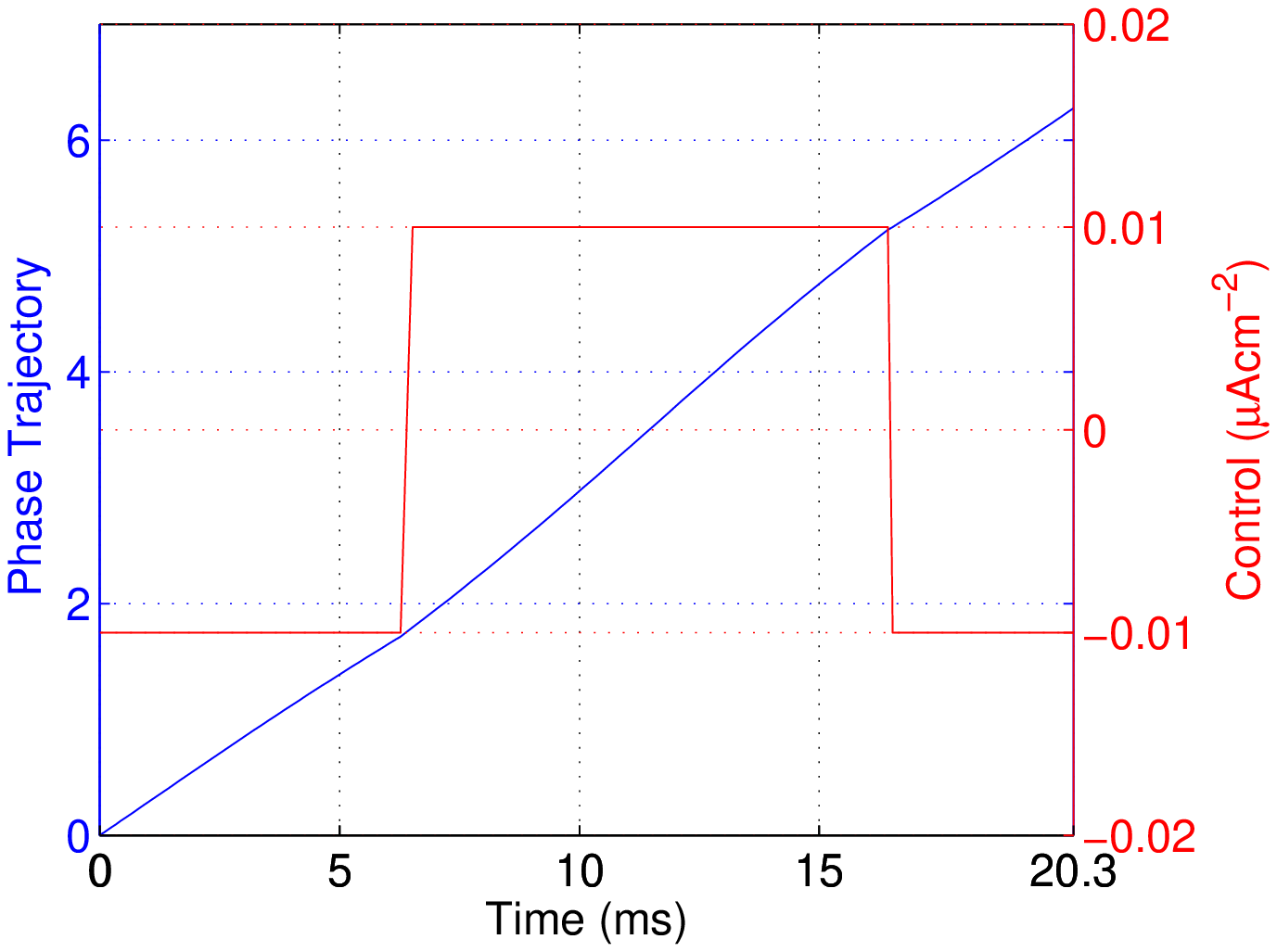} \label{fig:ML_min_time}} \\
		\hspace{0.7cm}
		\subfigure[]{\includegraphics[scale=0.5]{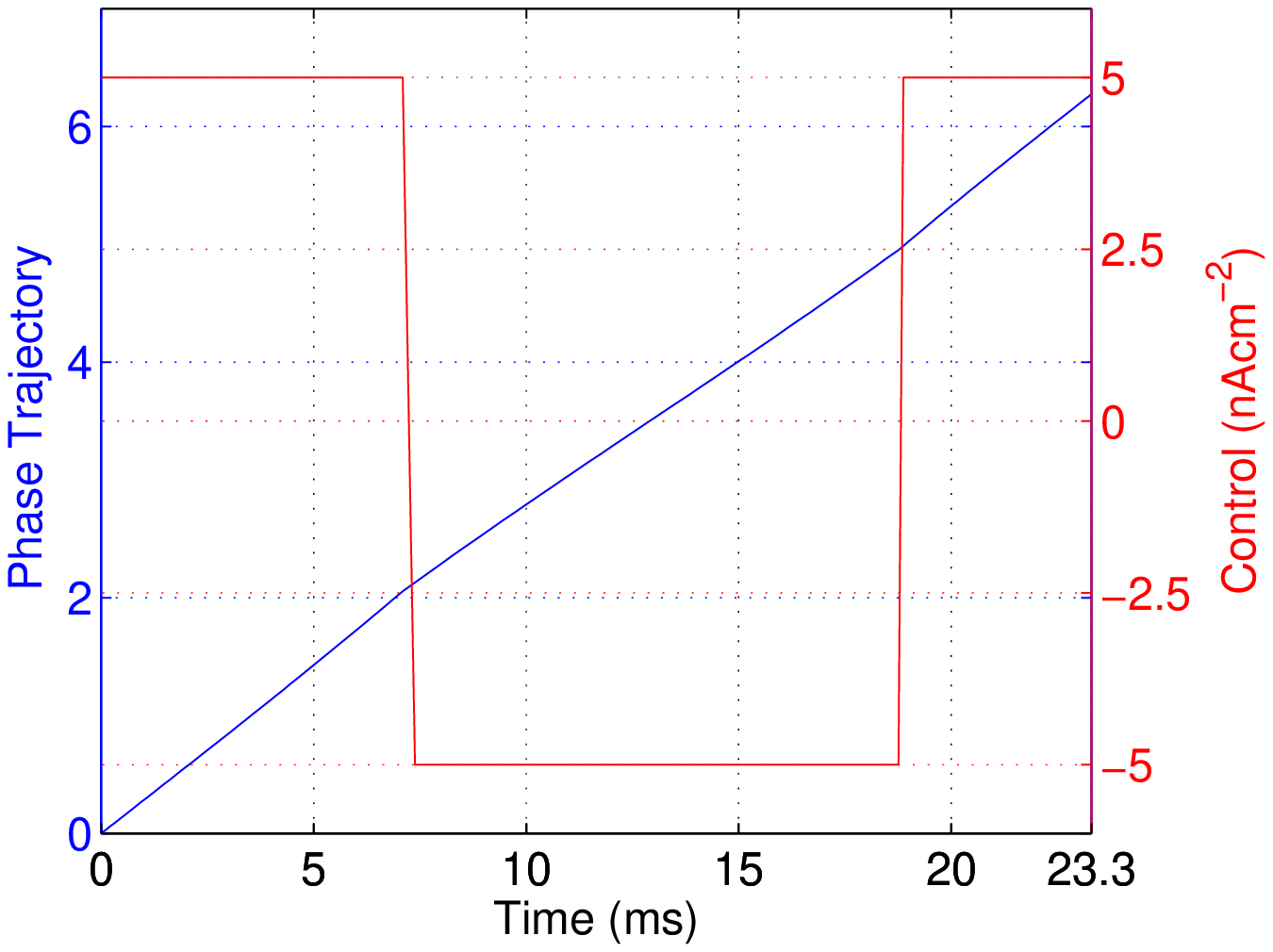}
		\label{fig:ML_max_time_bang_bang}} & 
		\subfigure[]{\includegraphics[scale=0.5]{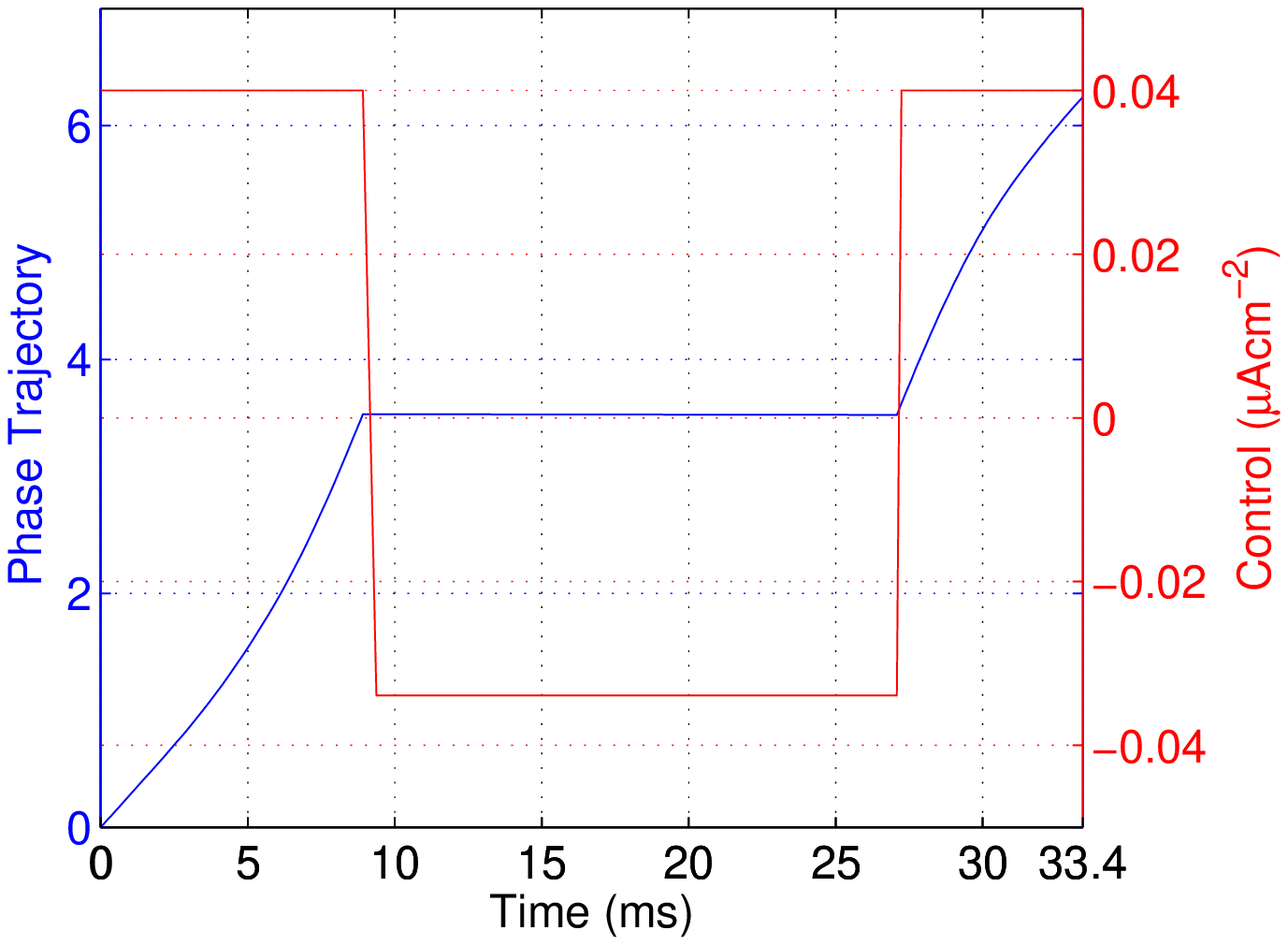}
		\label{fig:ML_max_time_bang_sing_bang}}
	\end{tabular}
	\caption{\ref{fig:ML_PRC} The Morris-Lecar PRC $Z(\theta)$ and its derivatives, $\frac{dZ}{d\theta}$ and $\frac{d^2Z}{d\theta^2}$. \ref{fig:ML_min_time} The charge-balanced minimum-time control and the corresponding phase trajectory for the Morris-Lecar phase model with respect to the bound on the control amplitude $M=0.01\ \mu Acm^{-2}$. \ref{fig:ML_max_time_bang_bang} and \ref{fig:ML_max_time_bang_sing_bang} show the charge-balanced maximum-time controls and the corresponding phase trajectories with $M=0.005\ \mu Acm^{-2}$ and  $M=0.04 \ \mu Acm^{-2}$, respectively.} 
\label{fig:ML_controls_and_trajectories}
\end{figure}

The Morris-Lecar neuron model was originally developed to capture the oscillatory behavior of barnacle muscle fibers \cite{Morris81}. It has been observed through experiments that the PRC for an Aplysia motoneuron is extremely similar to that of a Morris-Lecar PRC \cite{Milton00}. 
We consider a Morris-Lecar system with parameter values given in \cite{Dasanayake11}, which has a natural frequency $\omega=0.283\, rad/ms$. The PRC is approximated by \eref{eq:PRC_coeff} with the coefficients shown in \Tref{Table:ML_coefficient} and is illustrated, with its derivatives, in \Fref{fig:ML_PRC}.

\begin{table}[ht]
	\centering
		\begin{tabular}{|c|c|c|c|c|c|c|c|c|}
		\hline
		i & 1 & 2 & 3 & 4 & 5 & 6 & 7 & 8\\
		\hline
			$a_i$ & 5.137 & 5.773 &  0.7703 & 1.065 & 0.8143 & 0.1028 & 0.09711 & 0.0698 \\
			$b_i$ & 0.4356 & 0.7105 & 2.185 & 3.09 & 3.362& 4.876 & 5.829 & 6.525  \\
			$c_i$ & 1.005 & -1.474 & 0.6535 & 1.238 & 3.585 & 2.154  & 2.375 & 3.446 \\
			\hline
		\end{tabular}
		\caption{The coefficients of the equation \eref{eq:PRC_coeff} for the Morris-Lecar PRC}
		\label{Table:ML_coefficient}
\end{table}
Three examples are made to show the different structures of the optimal controls that are associated with various values of $M$ for the Morris-Lecar phase model. The charge-balanced minimum-time control and the resulting phase trajectory for $M=0.01\ \mu Acm^{-2}$ are given in \Fref{fig:ML_min_time}. The charge-balanced maximum-time controls and the respective trajectories subject to $M=0.005\ \mu Acm^{-2}$ and $M=0.04\ \mu Acm^{-2}$ are given in \Fref{fig:ML_max_time_bang_bang} and \Fref{fig:ML_max_time_bang_sing_bang}, respectively. The derivations of these optimal controls follow a similar procedure presented in \ref{Appendix_HH}.


\section{Validation of Phase Model Reduction to Full State-Space Model}
\label{sec:validation}
Because phase models of importance to applications are reductions of original higher dimensional state-space systems, we  explore in this section the extent to which controls synthesized using the former can achieve the desired objectives when applied to the latter. This will provide insight into the limits of the model reduction with respect to control synthesis, and allow the relationship to be calibrated for practical applications where the weak forcing assumption must be relaxed. Such an important validation is largely lacking in the literature. 

\begin{figure}
	\centering
	\includegraphics[scale=0.25]{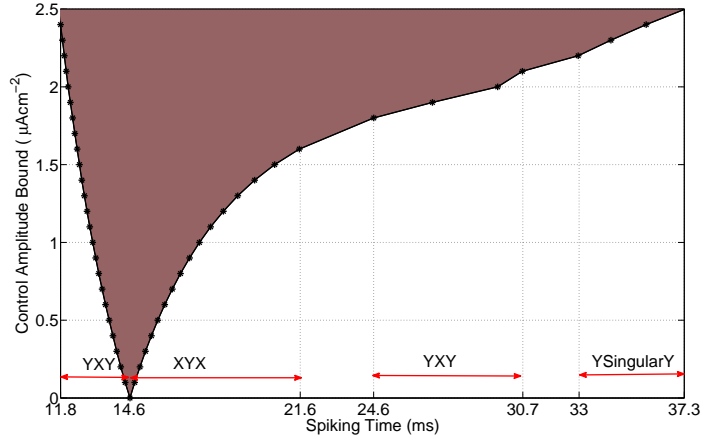}
	\caption{ A characterization of the realizable spiking times with respect to the bound on the control amplitude, $M\in[0,2.5]$, for the Hodgkin-Huxley phase model. The shaded region indicates the feasible spiking range resulting from the minimum- and maximum-time controls. Those minimum times (left to the natural spiking time $T_0=14.6\ ms$) are obtained by $YXY$ controls and maximum-times (right to $T_0$) are obtained by $XYX$, $YXY$ and $Y$-singular-$Y$ controls depending on $M$.}
	\label{fig:HH_feasible_spiking_time}
\end{figure}

We validate our optimal control strategies derived based on the phase models with the corresponding original state-space models. Specifically, we consider the Hodgkin-Huxley model. Note that an analytical derivation of the optimal controls directly from the state-space models is in general intractable and computationally expensive. A validation of the minimum and maximum spiking times with respect to the bound on the control amplitude is depicted in Figure \ref{fig:HH_feasible_spiking_time}, where the feasible spiking times are indicated as the shaded area.  
Each asterisk point on this graph represents the Hodgkin-Huxley neuron spiking time achieved by a particular form of the optimal control. The points correspond to minimum spiking times, which are less than the natural spiking time $T_0=14.64\ ms$, are obtained by $YXY$ controls, whereas the points correspond to maximum spiking times may be obtained by three structurally different controls, i.e., $XYX$, $YXY$, and $Y$-singular-$Y$ controls. 
This figure illustrates the limits on possible spiking times of the Hodgkin-Huxley model, which is important to the design of practical control inputs. For example, the knowledge of the feasible spiking range is helpful in designing optimal controls with other objectives such as minimum power controls \cite{Dasanayake11b}.

\begin{figure}[ht!]
	\centering
	\subfigure[Uncontrolled and controlled spiking trains for minimum time with amplitude $M=0.7\ \mu Acm^{-2}$ of Hodgkin-Huxley neuron.]{
\includegraphics[scale=0.3]{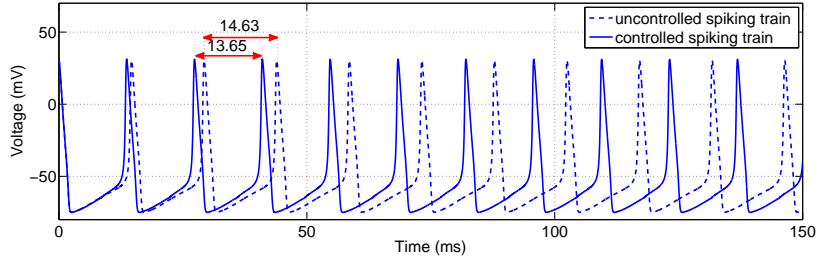}
	\label{fig:HH_simulation_min_time_M_07}
}
	\subfigure[Uncontrolled and controlled spiking trains for maximum time with  amplitude $M=0.7\ \mu Acm^{-2}$ of Hodgkin-Huxley neuron.]{
\includegraphics[scale=0.3]{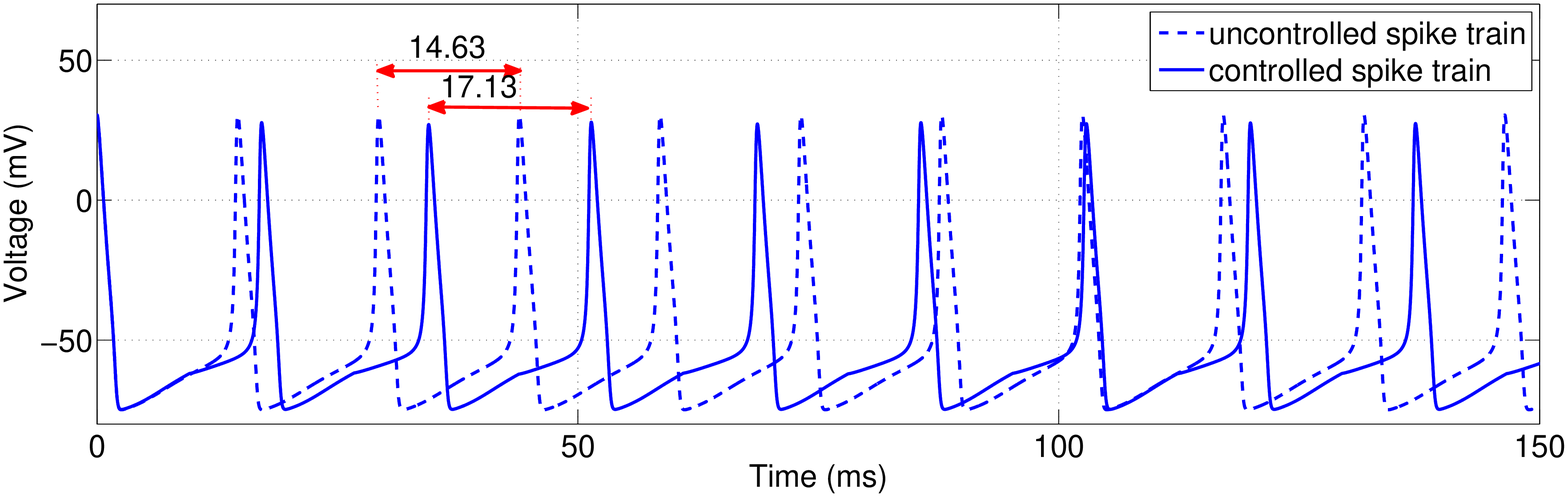}
	\label{fig:HH_simulation_max_time_M_07}
}
\caption{Application of derived optimal controls according to phase models to full Hodgkin-Huxley model}
\end{figure}

\begin{figure}
	\centering
		\includegraphics[scale=0.25]{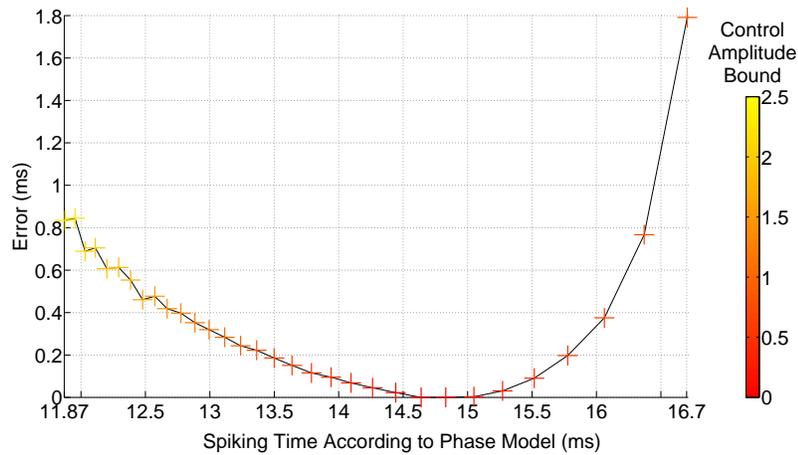}
		\caption{The absolute error in the spiking time when applying the charge-balanced time-optimal controls derived based on the Hodgkin-Huxley phase model to its full state-space model. The bound of the control amplitude is indicated as the color bar.}
	\label{fig:HH_error}
\end{figure}

The optimal controls derived based on the Hodgkin-Huxley phase model, shown in \Fref{fig:HH_min_time} and \ref{fig:HH_max_time_bang_bang}, are applied to the full Hodgkin-Huxley model, and a repeated application of such controls results in the desired spiking trains as displayed in \Fref{fig:HH_simulation_min_time_M_07} and \ref{fig:HH_simulation_max_time_M_07}. The respective minimum and maximum spiking times induced from these optimal controls subject to the amplitude bound $M=0.7\ \mu Acm^{-2}$ are $13.5\ ms$ and $16.37\ ms$ in the phase model and $13.65\ ms$ and $17.13\ ms$ in the full state-space model. Such an inconsistence is due to the model reduction, however, the resulting spiking behavior of the full Hodgkin-Huxley model shows great qualitative agreement with that of the phase model. The variation of the absolute errors between the actual and designed spiking times is shown in \Fref{fig:HH_error}, where the spiking behavior predicted based on the phase model matches better the full state-space model towards the weak forcing region. 

\section{Conclusion and Future Work}
\label{sec:future_work} 
In this paper, we investigated time-optimal controls for phase models of spiking neuron oscillators. In particular, we derived charge-balanced controls that lead to the minimum and the maximum inter-spike time of a neuron for a given bound on the control amplitude. 
 We showed that such optimal controls involve bang-bang and bang-singular-bang structures depending on the allowable control amplitude. 
 Although the amplitude level of weak forcing in the phase model is not practically quantifiable and can be greatly dependent on the dynamics of the system, our optimal control solutions were constructed for an arbitrary choice of bounds on the control amplitude, which accounts for this practical issue. We apply the derived optimal spike timing controls to  commonly-used phase models 
of oscillatory neurons to demonstrate their applicability to neuroscience. The methodology presented in this paper is general 
and can be applied not only to oscillatory neuron systems, but also to any oscillating system that can be represented by phase models including biological, chemical, electrical, and mechanical oscillators.

The theoretical results of this work characterize the fundamental limits on neuron spiking times that can be achieved by use of a charge-balanced bounded external input, and have potential impact on the improvement and development of innovative therapeutic procedures for neurological disorders. The extension of this work to the optimal control of a neuron population is of fundamental and practical importance. Our recent work has shown that an ensemble of uncoupled neurons is controllable and the minimum-power controls that spike a network of heterogeneous neurons can be constructed using a multidimensional pseudospectral method \cite{Li_NOLCOS10}. We plan to extend this current work to investigate the controllability and optimal controls of a network of coupled neurons.


\appendix
\section{ The Pontryagin's Maximum Principle}
\label{Appendix1}
\begin{thm}[Time-Optimal Control \cite{Pontryagin62}]
	\rm Let $(x_*(t), u_*(t))$ be a time-optimal controlled trajectory that transfers the initial
condition $x(0) = x_0$ into the terminal state $x(T ) = x_T$ . Then, it is a necessary condition for optimality that there exists a constant $\lambda_0\geq0$ and nonzero, absolutely continuous row vector function $\lambda(t)$ such that:
	\begin{enumerate}
		\item $\lambda$ satisfies the so-called adjoint equation
			\begin{equation*}
				\dot{\lambda}(t)=-\frac{\partial H}{\partial x}(\lambda_0,\lambda(t),x_*(t),u_*(t))
			\end{equation*}
		\item For $0 \leq t \leq T$ the function $u\mapsto H(\lambda_0, \lambda(t), x_*(t), u)$ attains its minimum over the control set $U$ at $u = u_*(t)$.
		\item $H(\lambda_0, \lambda(t), x_*(t), u_*(t)) \equiv 0$.
	\end{enumerate}
\end{thm}


 \section{ The Derivation of Time-Optimal Controls for the Hodgkin-Huxley Phase Model}
 \label{Appendix_HH}


\subsection*{Charge-Balanced Minimum-Time Control for Hodgkin-Huxley Phase Model}
 The Hodgkin-Huxley PRC given in Figure \ref{fig:HH_PRC} has at most two singular trajectories (points), $\theta_{s,1}=3.34$ and $\theta_{s,2}=4.58$,  calculated by the condition $\frac{\partial Z(\theta)}{\partial\theta}=0$. According to the shape of this PRC, there exist at most two switching points satisfying $Z(\t)=\alpha$, where $\alpha$ is a constant defined in \eref{eq:switch_modified}. Since the minimum-time control takes the bang-bang form as shown in Section \ref{sec:deriv_min}, it requires to calculate the switching points and determine the type of the switching at these points for the optimal control synthesis. At the switching points, $\dot{\phi}=-\partial Z/\partial\t$ is given by \eref{eq:dot_phi_at_switch}, and hence a $Y$ to $X$ switch may occur in the region $R_1=[0,\theta_{s,1}]$ or $R_3=[\theta_{s,2},2\pi]$, and an $X$ to $Y$ switch may occur in $R_2=[\theta_{s,1},\theta_{s,2}]$. This implies that bang-bang controls with one switch, such as the $XY$ or $YX$ form, are not feasible solutions because these controls will violate the charge-balance constraint. Consequently, the optimal control has two switching points, and the candidate is either a $YXY$ trajectory with one switch in the interval $R_1$ and one in $R_2$, or an $XYX$ trajectory with one switch in $R_2$ and one in $R_3$. We can further simplify the possible intervals of switching by observing the shape of the PRC. The Hodgkin-Huxley PRC depicted in \Fref{fig:HH_PRC} has three zeros at $\theta_{r,1}=0$, $\theta_{r,2}=3.86$, and $\theta_{r,3}=2\pi$. Therefore, for an optimal $YXY$ trajectory the first and the second switch will occur in $[0,\theta_{s,1}]$ and  $[\theta_{s,1}, \theta_{r,2}]$, respectively, and for an optimal $XYX$ trajectory, they will occur in $[\theta_{r,2},\theta_{s,2}]$ and $[\theta_{s,2},2\pi]$, respectively. The minimum-time control is then selected between these two. Note that for a given bound $M$, it may not be possible to have both $XYX$ and $YXY$ solutions. For example, if the bound is $M=0.7$, then the only feasible optimal solution is $YXY$. In this case, the two switching points $\theta_1$ and $\theta_2$ can be calculated through

\begin{eqnarray}
	\label{eq:sp1}
	&0=\int_0^{\theta_1}{\frac{M}{\omega+MZ(\theta)}}d\theta+\int_{\theta_1}^{\theta_2}{\frac{-M}{\omega-MZ(\theta)}}d\theta+\int_{\theta_2}^{2\pi}{\frac{M}{\omega+MZ(\theta)}}d\theta,\\
	\label{eq:sp2}
	&Z(\theta_1)=Z(\theta_2),
\end{eqnarray}
and the control is then given by
\begin{equation*}
\label{eq:min_time_control_HH_singular}
    u^*_{min}=\left\{\begin{array}{ll} M, & 0\leq\theta\leq\theta_1 ,\\ -M, & \theta_1<\theta<\theta_2, \\ M, & \theta_2\leq\theta\leq2\pi.
    \end{array}\right.
\end{equation*}

\subsection*{Charge-Balanced Maximum-Time Control for Hodgkin-Huxley Phase Model}
In the case of the maximum-time control, the two singular points, $\theta_{s,1}$ and $\theta_{s,2}$, are candidates for the optimal trajectory because they are slower than the bang trajectories as proved in Section \ref{sec:max_time_sniper}. Letting $\dot{\theta}=0$ in \eref{eq:phasemodel}, we find the controls that keep the trajectory at the  singular points are $u_{s,1}=-\frac{\omega}{Z(\theta_{s,1})}=3.50$ and $u_{s,2}=-\frac{\omega}{Z(\theta_{s,2})}=-2.15$. 
 There exist three cases when constructing maximum-time controls according to $M$ and thus to the feasibility of $u_{s,1}$ and $u_{s,2}$.

\emph{(Case I: $M<|u_{s,2}|$)}
In this case, both the singular points $\theta_{s,1}$ and $\theta_{s,2}$ are infeasible. Therefore, the optimal control is bang-bang and is in fact the opposite of the minimum-time control described above. Similar to the minimum-time case, we can calculate the corresponding $XYX$ and $YXY$ solutions and choose the maximum time achieved between these scenarios. For example, consider the bound $M=0.7$, then the only solution is $XYX$ and the two switching points are calculated by substituting $M$ with $-M$ in \eref{eq:sp1} and solving \eref{eq:sp1} and \eref{eq:sp2}.
The optimal bang-bang control is then given by
\begin{equation*}
\label{eq:max_time_control_HH}
    u^*_{max}=\left\{\begin{array}{ll} -M, & 0\leq\theta<\theta_1 ,\\ M, & \theta_1\leq\theta\leq\theta_2, \\ -M, & \theta_2<\theta\leq2\pi.
    \end{array}\right.
\end{equation*}
\emph{(Case II: $|u_{s,2}|\leq M <|u_{s,1}|$)}
In this case, $\theta_{s,2}$ is the only feasible singular trajectory (point) generated by the singular control $u_{s,2}=-2.15<0$.
Because there are only two switching points allowed in the optimal trajectory, this together with the fact that $u_{s,2}$ is of negative charge forces the optimal control to take the ``$Y$-singular-$Y$'' form given by
\begin{equation*}
	\label{eq:max_time_control_HH_singular}
	u^*_{max}=\left\{\begin{array}{ll} M, & 0\leq\theta< \theta_{s,2} ,\\ u_{s,2}, & \theta=\theta_{s,2}, \\ M, & \theta_{s,2}<\theta\leq2\pi.\end{array}\right.
\end{equation*}
Similar to the SNIPER phase model described in Section \ref{sec:max_time_sniper}, the time it takes to reach the singular point is given by,
\begin{equation*}
	t_1=\int_0^{\theta_{s,2}}{\frac{1}{\omega+Z(\theta)M}}d\theta
\end{equation*}
and the time required to reach the target point $2\pi$ from the point $\theta_{s,2}$ is
\begin{equation*}
	t_3=\int_{\theta_{s,2}}^{2\pi}{\frac{1}{\omega+Z(\theta)M}}d\theta.
\end{equation*}
The time during which the trajectory stays on $\t_{s,2}$ is determined by the charge-balance constraint 
and is given by 
\begin{equation*}
	t_2=\left|\frac{(t_1+t_3)M}{u_{s,2}}\right|.
\end{equation*}
Now, the optimal control can be stated in terms of time as
\begin{equation}
	\label{eq:max_time_control_HH_singular_vs_time}
    u^*_{max}=\left\{\begin{array}{ll} M, & 0\leq t< t_1 ,\\ u_{s,2}, & t_1\leq t\leq t_1+t_2, \\ M, & t_1+t_2<t\leq t_1+t_2+t_3.
\end{array}\right.
\end{equation}

\emph{(Case III: $|u_{s,1}|\leq M$)}
In this case, staying on either singular point is possible by using an appropriate control. Furthermore, since the two singular controls have opposite signs, the charge-balance constraint can be preserved by staying for an appropriate time period at each singular point. As a result, the spiking time can be arbitrarily delayed, which may not be of practical interest due to the requirement of  relatively high amplitude.

\section*{References}
\bibliographystyle{vancouver}	
\bibliography{cbtoc}

\begin{thebibliography}{10}

\bibitem{Izhikevich07}
Izhikevich IE.
\newblock Dynamical Systems in Neuroscience.
\newblock Cambridge, Massachusetts: The MIT Press; 2007.

\bibitem{Schiff94}
Schiff SJ, K J, H DD, T C, L SM, L WD.
\newblock Controlling chaos in the brain.
\newblock Nature. 1994;370:615--620.

\bibitem{Benabid10}
Benabid AL, Pollak P, Haffmann D, Gervason C, Hommek M, Perret JE, et~al.
\newblock Long-term suppression of tremor by chronic stimulation of the ventral
  intermediate thalamic nucleus.
\newblock The Lancet. 1991;337:403--406.

\bibitem{Lozano04}
Lozano AM, Eltahawy H.
\newblock How does DBS work?
\newblock Suppl Clin Neurophysiol. 2004;57:733--736.

\bibitem{Nabi11}
Nabi A, Moehlis J.
\newblock Single input optimal control for globally coupled neuron networks.
\newblock J Neural Eng. 2011;8:065008.

\bibitem{Ortmanns07}
Ortmanns M.
\newblock Charge Balancing in Functional Electrical Stimulation: A Comparative
  Study.
\newblock In: Proc. of IEEE International Symposium on Circuits and Systems.
  IEEE; 2007. .

\bibitem{Merrill05}
Merril DR, Bikson M, Jefferys JGR.
\newblock Electrical stimulation of excitable tissue: design of efficacious and
  safe protocols.
\newblock Journal of Neuroscience Methods. 2005;141:717--198.

\bibitem{Brown04}
Brown E, Moehlis J, Holmes P.
\newblock On the Phase Reduction and Response Dynamics of Neural Oscillator
  Populations.
\newblock Neural Comput. 2004;16(4):673--715.

\bibitem{Ashwin92}
Ashwin P, Swift J.
\newblock The dynamics of n weakly coupled identical oscillators.
\newblock J Nonlinear Sci. 1992;2:69--108.

\bibitem{Tass89}
Tass PA.
\newblock Phase Resetting in Medicine and Biology.
\newblock New York: Springer; 1989.

\bibitem{Moehlis06}
Moehlis J, Shea-Brown E, Rabitz H.
\newblock Optimal Inputs for Phase Models of Spiking Neurons.
\newblock J Comput Nonlinear Dynam. 2006;1(4):358--367.

\bibitem{Dasanayake11}
Dasanayake I, Li JS.
\newblock Optimal design of minimum-power stimuli for phase models of neuron
  oscillators.
\newblock Phy Rev E. 2011;83:061916.

\bibitem{Dasanayake11b}
Dasanayake I, Li JS. Charge-Balanced Minimum-Power Controls for Spiking Neuron
  Oscillators; 2012.
\newblock arXiv:1111.2879.

\bibitem{Zlotnik12}
Zlotnik A, Li JS.
\newblock Optimal Entrainment of Neural Oscillator Ensembles.
\newblock J Neural Eng. 2012;9:046015.

\bibitem{Zlotnik11}
Zlotnik A, Li JS.
\newblock Optimal asymptotic entrainment of phase-reduced oscillators.
\newblock In: Proc. ASME Dynamic Systems and Control Conf. (Arlington, VA).
  ASME; 2011. p. 479--84.

\bibitem{kiss10}
Harada T, Tanaka HA, Hankins MJ, Kiss IZ.
\newblock Optimal Waveform for the Entrainment of a Weakly Forced Oscillator.
\newblock Phys Rev Lett. 2010 Aug;105(8):088301.

\bibitem{Stefanatos12}
Stefanatos D, Li JS.
\newblock Antiphase synchronization of phase-reduced oscillators using
  open-loop control.
\newblock Phys Rev E. 2012 Mar;85:037201.

\bibitem{Feng03}
Feng J, Tuckwell HC.
\newblock Optimal Control of Neuronal Activity.
\newblock Phy Rev Lett. 2003;91(1):018101.

\bibitem{Li_NOLCOS10}
Li JS.
\newblock Control of a Network of Spiking Neuron.
\newblock In: Proc. 8th IFAC Symposium on Nonlinear Control Systems, Bologna,
  Italy, 2010. IFAC; 2010. .

\bibitem{Li12}
Li JS, Dasanayake I, Ruths J. Control and Synchronization of Neuron Ensembles;
  2012.
\newblock arXiv:1111.6306v2.

\bibitem{Nabi10}
Nabi A, Moehlis J.
\newblock Time Optimal Control of Spiking Neurons.
\newblock J Math Biol. 2011;7(6):2466--2470.

\bibitem{Danzl08}
Danzl P, Moehlis J.
\newblock Spike Timing Control of Oscillatory Neuron Models Using Impulsive and
  Quasi-Impulsive Charge-Balanced Inputs.
\newblock In: Proc. American Control Conference, Seattle, Washington, USA,
  2008. IEEE; 2008. .

\bibitem{Hodgkin52}
Hodgkin AL, Huxley AF.
\newblock A Quantitative Description of Membrane Current and Its Application To
  Conduction and Excitation in Never.
\newblock J Physiol. 1952;117:500--544.

\bibitem{Pontryagin62}
Pontryagin LS, Boltyanskii VG, Gamkrelidze RV, Mishchenko EF.
\newblock The Mathematical Theory of Optimal Processes.
\newblock John Wiley \& Sons, Inc.; 1962.

\bibitem{Morris81}
Morris C, Lecar H.
\newblock Voltage oscillations in the barnacle giant muscle fiber.
\newblock Biophys J. 1981;35(1):193--213.

\bibitem{Bonnard03}
Bonnard B, Chyba M.
\newblock Singular Trajectory and their Role in Control Theory.
\newblock New York: Springer-Verlag; 2003.

\bibitem{Ermentrout96}
Ermentrout B.
\newblock Type {I} membranes, phase resetting curves, and synchrony.
\newblock Neural Comput. 1996;8:979--1001.

\bibitem{Rose89}
Rose R, Haindmarsh J.
\newblock The assembly of ionic currents in a thalamic neuron I. The
  three-dimensional model.
\newblock Proc R Soc Lond B. 1989;237:267--288.

\bibitem{Milton00}
Foss J, Milton J.
\newblock Multistability in Recurrent Neural Loops Arising From Delay.
\newblock J Neurophysiol. 2000;84:975--985.

\end{thebibliography}
\end{document}